\documentclass[a4paper,11pt,reqno,noindent]{amsart}
\usepackage[centertags]{amsmath}
\usepackage{amsfonts,amssymb,amsthm} 
\usepackage[english]{babel}
\usepackage{newlfont}
\usepackage{color}
\usepackage[body={15cm,21.5cm},centering]{geometry} 
\usepackage{fancyhdr}
\usepackage{esint}

\usepackage{enumerate}

\newtheorem{theo}{Theorem}

\newtheorem{prop}{Proposition}
\newtheorem{coro}{Corollary}

\theoremstyle{definition}
\newtheorem{rem}{Remark}

\newtheorem{hypo}{Hypothesis}

\numberwithin{equation}{section}
\numberwithin{lemma}{section} 
\numberwithin{defi}{section} 
\numberwithin{rem}{section} 

\newcommand{\per}{\mathrm{per}}

\newcommand{\R}{\mathbb{R}}
\newcommand{\Z}{\mathbb{Z}}
\newcommand{\Q}{\mathbb{Q}}
\newcommand{\N}{\mathbb{N}}

\newcommand{\calF}{\mathcal{F}}

\newcommand{\calE}{\mathcal{E}}

\newcommand{\calL}{\mathcal{L}}

\mathchardef\emptyset="001F

\newcommand{\var}[1]{\mathrm{var}\left[#1\right]}
\newcommand{\cov}[2]{\mathrm{cov}\left[#1;#2\right]}
\newcommand{\expec}[1]{\left\langle #1 \right\rangle}

\newcommand{\step}[1]{\noindent \textit{Step} #1.}

\title[Stationary random point sets]
{When are increment-stationary random point sets stationary?}
\author[A. Gloria]{Antoine Gloria}
\date{\today}
\address[Antoine Gloria]{Universit\'e Libre de Bruxelles (ULB) \\ Brussels, Belgium \\ and Team MEPHYSTO \\  Inria Lille - Nord Europe \\ Villeneuve d'Ascq, France}
\email{agloria@ulb.ac.be}
\begin{document}
\maketitle

\begin{center}
\begin{minipage}{13cm}
\small{
\noindent {\bf Abstract.} 
In a recent work, Blanc, Le Bris, and Lions defined a notion of increment-stationarity for random point sets, which allowed them to prove the existence of a thermodynamic limit for two-body potential energies on such point sets (under the additional assumption of ergodicity), and 
to introduce a variant of stochastic homogenization for increment-stationary coefficients.
Whereas stationary random point sets are increment-stationary, it is not clear a priori under which conditions increment-stationary random point sets are stationary.
In the present contribution, we give a characterization of the equivalence of both notions of stationarity based on elementary PDE theory in the probability space.
This allows us to give  conditions on the decay of a covariance function associated with the random point set, which ensure that increment-stationary random point sets are stationary random point sets up to a random translation with bounded second moment in dimensions $d>2$. In dimensions $d=1$ and $d=2$, we show that such sufficient conditions cannot exist.
\vspace{10pt}


\vspace{6pt}
}

\end{minipage}
\end{center}

\bigskip


\section{Introduction and main results}

\subsection{Increment-stationarity}\label{sec:1d}

Let $d$ be the dimension, and denote by $\calL(\R^d)$ the set of locally finite simple point sets of $\R^d$.
In what follows we consider random variables on some probability space $(\Omega,\mathcal{F},\mathbb{P})$
taking values in $\calL(\R^d)$. We call such random variables random point sets.
In \cite{BLL-07}, Blanc, Le Bris and Lions addressed the issue of defining the thermodynamic limit of the energy of random sets $\ell$ of particles (seen as simple random point sets).
Typical energies to be considered are given by two-body potentials $V:\R^d\setminus\{0\}\to \R$, so that the energy in some bounded region $D\subset \R^d$ writes
$$
\calE(\ell,D)\,:=\,\frac{1}{2}\sum_{x,y\in \ell\cap D, x\neq y}V(x-y).
$$
As noticed in \cite{BLL-07,BLL-07b}, the existence of a deterministic thermodynamic limit
\begin{equation}\label{eq:therm-limit}
\lim_{D\uparrow \R^d} \frac{1}{|D|} \calE(\ell,D)
\end{equation}
(when $D$ properly invades $\R^d$) does not require the point set $\ell$ to be stationary, which lead the authors to define a notion of increment-stationarity.

\medskip

Let us start by recalling the definition of increment-stationarity.
We say that $\ell$ is \textit{increment-stationary} if there exist a sequence of random vectors $\{Y_k\}_{k\in \Z^d}$ of $L^2(\Omega,\R^d)$
and a measure-preserving group action $\{\theta_k\}_{k\in \Z^d}$ such that $\ell$ satisfies for almost all $\omega\in \Omega$ and all $k\in \Z^d$,
\begin{equation}\label{eq:def-incr-stat-1}
\ell(\theta_k \omega)\,=\,\ell(\omega)+Y_k(\omega)
\end{equation}
that is, $\ell(\theta_k \omega)$ is the translation of $\ell(\omega)$ by the vector $Y_k(\omega)$. 
Assume in addition that
\begin{eqnarray}
Y_0&\equiv&0,\label{eq:def-incr-stat-2}\\
Y_{j}(\theta_k\omega)-Y_{i}(\theta_k\omega)&=&Y_{j+k}(\omega)-Y_{i+k}(\omega)\label{eq:def-incr-stat-3}
\end{eqnarray}
for all $i,j,k\in \Z$ and almost every $\omega\in \Omega$ (properties which directly follow from \eqref{eq:def-incr-stat-1} if
$Y_k$ is uniquely defined by \eqref{eq:def-incr-stat-1}).
Let us check that such a definition ensures that the limit in \eqref{eq:therm-limit} is deterministic if it exists.
For almost all $\omega$, $Y_k(\omega)$ is finite for all $k\in \Z^d$.
Hence, for all $D\subset \R^d$, the points of $\ell(\theta_k \omega)\cap D$ are the translations of the points 
of $\ell(\omega)\cap (D-Y_k(\omega))$. Since the two-body potential is non-negative and translation-invariant, this yields
\begin{equation*}
 \calE(\ell(\theta_k\omega),D) = \calE(\ell(\omega),D-Y_k(\omega)) \,\leq \, \calE(\ell(\omega),B_{|Y_k(\omega)|}(D)),
\end{equation*}
where for all $t>0$, $B_t(D):=\{x\in \R^d, d(x,\partial D)<t\}$.
Likewise,
\begin{equation*}
 \calE(\ell(\omega),D) = \calE(\ell(\theta_k\omega),D+Y_k(\omega)) \,\leq \, \calE(\ell(\theta_k \omega),B_{|Y_k(\omega)|}(D)).
\end{equation*}
Hence, since for all $t\ge 0$,  $\lim_{D\uparrow \R^d}\frac{|B_t(D)|}{|D|}=1$, these two estimates yield for all $k\in \Z^d$,
$$
\limsup_{D\uparrow \R^d} \frac{1}{|D|} \calE(\ell(\omega),D)\,=\, \limsup_{D\uparrow \R^d} \frac{1}{|D|} \calE(\ell(\theta_k\omega),D),
$$
and the thermodynamic limit is invariant by the group action, and therefore constant if the group action is ergodic.

We say that $\ell$ is \textit{stationary} if it is increment-stationary and if $Y_k$ is given for all $k\in \Z^d$ by
\begin{equation}
Y_k:=Tk
\end{equation}
for some deterministic $d\times d$ matrix $T$.

In \cite{BLL-07}, the authors prove that if the positions of an infinite set of nuclei are given by a stationary random point set (satisfying in addition uniform hard-core and non-empty space properties), then the thermodynamic limit of 
the associated electronic cloud exists in the sense that the notions of averaged energy and cloud density are well-defined, in the case
of Thomas-Fermi models.  
This was later extended by Blanc and Lewin \cite{Blanc-Lewin-12}, and Canc\`es, Lahbabi and Lewin \cite{CLL-13}, to quantum models with Coulomb forces and to Hartree-Fock and Kohn-Sham type models, respectively.
In terms of point sets, their proofs essentially rely on the stationarity of two-body interactions and an ergodic theorem, so that, using the argument above, these proofs should extend to the more general case of increment-stationary random point sets.

\medskip

The aim of this contribution is to investigate in which respect \textit{increment-stationarity} is more general than \textit{stationarity}, and identify under which conditions on the sequence $Y_k$ one can conclude that 
an \textit{increment-stationary} random point set is \textit{stationary}.
In particular, both the probability space and the group action are fixed, and we are indeed investigating the \textit{rigidity} of increment-stationary random point sets.
As we shall see, the validity of such a rigidity result depends on the dimension.

\subsection{Main results}

In what follows we endow $(\Omega,\mathcal{F},\mathbb{P})$ with 
an ergodic measure-preserving discrete group action $\{\theta_k\}_{k\in \Z^d}:\Omega\to \Omega$, 
and we denote by $\expec{\cdot}$, $\var{\cdot}$, and $\cov{\cdot}{\cdot}$ the associated expectation, variance, and covariance, respectively.
We denote by $L^2_0(\Omega,\R^d)$ the space of random vectors with bounded second moments and vanishing expectations.
We let $\{e_l\}_{l\in \{1,\dots,d\}}$ denote the canonical basis of $\R^d$.

\medskip

Let $\ell$ be a stationary point set in $\R^d$ and $Y\in L^2_0(\Omega,\R^d)$ be some non-identically zero random vector.
Then $\ell+Y$ is not stationary but is clearly increment-stationary. 
We shall say that a random point set $\ell$ is \textit{stationary up to translation} if there exists some $\tilde Y\in L^2_0(\Omega,\R^d)$ 
such that $\ell+\tilde Y$ is a stationary random point set.
The following 
theorem gives a characterization of the equivalence between \textit{increment-stationarity} and \textit{stationarity up to translation}.
This result is directly inspired by the treatment of the corrector equation in stochastic homogenization by Papanicolaou and Varadhan in \cite{Papanicolaou-Varadhan-79} (see \cite{Kunnemann-83} for the case of discrete elliptic equations).
It relies on the differential calculus in the probability space generated by the group action $\{\theta_k\}_{\theta\in \Z^d}$.
\begin{prop}\label{prop}
Let $\ell:\Omega\to \calL(\R^d)$ be an increment-stationary random point set for the group action $\{\theta_k\}_{k\in \Z^d}$. Let $\{Y_k\}_{k\in \Z^d} \in L^2(\Omega,\R^d)$ be a sequence satisfying \eqref{eq:def-incr-stat-1}---\eqref{eq:def-incr-stat-3}. For all $\mu>0$ and $i\in \{1,\dots,d\}$, consider the unique weak solution $\phi_{\mu,i} \in L^2(\Omega)$ of the equation: For all $\psi\in L^2(\Omega)$,
\begin{equation}\label{eq:corrector}
\expec{\mu \phi_{\mu,i}\psi+D\psi\cdot D\phi_{\mu,i}}\,=\,\expec{D\psi\cdot \zeta_i},
\end{equation}
where $\zeta_i \,:=\, (Z_{e_1}\cdot e_i,\dots,Z_{e_d}\cdot e_i)\in L^2(\Omega,\R^d)$,
$Z_{e_l}:=Y_{e_l}-\expec{Y_{e_l}}$ for all $l\in \{1,\dots,d\}$, and $D:=(D_1,\dots,D_d)$ is the differential operator from $L^2(\Omega)\to L^2(\Omega,\R^d)$ defined by $D_l \psi(\omega)=\psi(\theta_{e_l}\omega)-\psi(\omega)$.
If $\{\phi_{\mu,i}\}_{i\in \{1,\dots,d\}}$ is bounded in $L^2(\Omega)$ uniformly wrt $\mu>0$, then $\ell$ is stationary up to translation.
Conversely, if $\ell$ is stationary up to translation and if $\ell$ is non-degenerate (in the sense for almost all $\omega$, if $X$ is such that  $\ell(\omega)+X=\ell(\omega)$ then $X=0$), then 
$\{\phi_{\mu,i}\}_{i\in \{1,\dots,d\}}$ is bounded in $L^2(\Omega)$ uniformly wrt $\mu>0$.
\end{prop}
\begin{rem}
The condition that $\ell$ be non-degenerate ensures that the sequence $Y_k$ is unique and rules out periodicity.
In this case, \eqref{eq:def-incr-stat-2} and \eqref{eq:def-incr-stat-3} follow from \eqref{eq:def-incr-stat-1}.
If the point set is degenerate, the sequence $Y_k$ is not uniquely defined due to some translation-invariance.
\end{rem}
It remains to identify sufficient conditions on $\{Y_{e_l}\}_{l\in \{1,\dots,d\}}$ for the boundedness of the functions $\phi_{\mu,i}$. 
These conditions are written in terms of the decay of a covariance function as follows.
\begin{hypo}[Decay of order $\alpha>0$]\label{hypo}
The random point set $\ell:\Omega\to \calL(\R^d)$ is increment-stationary for the group action $\{\theta_k\}_{k\in \Z^d}$, and
the associated random vectors $Y_{e_1},\dots,Y_{e_d}$ display the following covariance decays: There exists $\alpha>0$ such that for all $k \in \Z^d$ and $l,l',n\in \{1,\dots,d\}$
$$
\cov{Y_{e_l}\circ\theta_k \cdot e_n}{Y_{e_{l'}}\cdot e_n} \,\lesssim \,\frac{1}{1+|k|^\alpha},
$$
where $\lesssim$ means $\leq$ up to a universal multiplicative constant.
\end{hypo}
As the following result shows, there are two types of behavior, depending on the dimension $d$ ($d=2$ is critical).
\begin{theo}\label{th:stat} We have:
\begin{itemize}
\item For $d=1$ and $d=2$: There exists an increment-stationary random point set $\ell:\Omega\to \calL(\R^d)$ satisfying Hypothesis~\ref{hypo} with a finite range of dependence and which is not stationary up to translation.
\item For all $d>2$: If $\ell:\Omega\to \calL(\R^d)$ is increment-stationary  and satisfies Hypothesis~\ref{hypo} for some $\alpha>2$, then $\ell$ is stationary up to translation.
\end{itemize}
\end{theo}
Let us comment on this result. On the one hand, there is a rather natural similarity with the behavior of the corrector in stochastic homogenization. On the other hand, there is also some connection with the Palm-Khinchin theory for point processes on the real line. 

We start with the connection to stochastic homogenization. In this case the gradient of the corrector is stationary (which corresponds to the increment-stationarity of the point set), and we investigate whether the corrector can be stationary itself. In dimension $d=1$, there cannot exist stationary correctors in $L^2(\Omega)$ (for this would contradict the central limit theorem, see discussion \cite[p.~790]{Gloria-Otto-09}), dimension $d=2$ is critical 
and stationary correctors do not exist either, whereas in dimensions $d>2$, stationary correctors exist under some assumptions on the statistics (a spectral gap estimate, see \cite{Gloria-Otto-09}).
Equation~\eqref{eq:corrector} can indeed be seen as the corrector equation in the regime of vanishing ellipticity contrast (the variable-coefficients elliptic operator is replaced by a constant-coefficients elliptic operator).
In particular, in dimensions $d=1$ and $d=2$, the corrector in the regime vanishing ellipticity contrast for
independent and identically distributed random conductivities provides with an example of increment-stationary point set which is not stationary up to translation,
see Step~2 in the proof of Theorem~\ref{th:stat} for details.

Let us now turn to the connection to the Palm-Khinchin theory. Let $d=1$, and consider a random point set $\ell$ on the real line. Recall (see for instance \cite[Chapter~3]{Daley-Vere-03}) that the random point set $\ell$ is characterized by
the sequence of (measurable) random variables $\{X_i\}_{i\in\Z}$, where $X_i<X_{i+1}$ and $X_0$ is the closest point to the origin on the negative axis. The associated sequence of intervals is denoted by $\{\tau_i\}_{i\in \Z^d}$ and defined by $\tau_i=X_{i+1}-X_i$.
The random point set $\ell$ is said to be \textit{interval-stationary} if the sequence $\{\tau_i\}_{i\in \Z}$ is stationary in the following sense: For all $m\in \N$, $i_1,\dots,i_m\in \Z$, the distribution of $(\tau_{i_1+k},\dots,\tau_{i_m+k})$ does not depend on $k\in \Z$.
In turn, this implies (and is indeed equivalent to, up to changing the probability space, see for instance such a construction in \cite[Section~16.1]{Koralov-Sinai-07}) the existence of a discrete group action $\{\theta_z\}_{z\in \Z}$ which preserves the probability measure and is such that for all $i,k\in \Z$ and almost every $\omega \in \Omega$,
$$
\tau_i(\theta_k \omega)\,=\,\tau_{i+k}(\omega).
$$
In terms of the random point set $\ell$, this turns into: For all $i,j,k\in \Z$ and almost every $\omega \in \Omega$,
$$
X_{i+k}(\omega)-X_{j+k}(\omega)\,=\,X_i(\theta_k\omega)-X_j(\theta_k \omega).
$$
The latter implies that the random point set $\ell$ is increment-stationary with $Y_k(\omega)=X_0(\theta_k\omega)-X_0(\omega)$. The Palm-Khinchin theory (see for instance \cite[Theorem~13.3.I]{Daley-Vere-03}) establishes a one-to-one relation between interval-stationary and stationary random point sets, the proof of which shows that to pass from one to the other one has to consider some random translation which \emph{is not} in $L^2(\Omega)$. The argument behind this is best illustrated by the ``waiting time paradox" (see e.~g. \cite[p6]{Daley-Vere-03}), which shows that the Poisson point process, which is obviously stationary, is not interval-stationary (the probability that 0 belongs to a large interval is larger than the one it belongs to a small interval, so that the length of the interval around the origin cannot be exponentially distributed: this interval has to be sent to infinity to turn the Poisson process into an interval-stationary process). The general incompatibility between stationary and interval-stationary point sets gives another interpretation of the example of Theorem~\ref{th:stat} for $d=1$. Define $\ell=\{X_k\}_{k\in \Z}$  by 
$X_0=0$, and for all $k\in \N$ by $X_{k}=X_{k-1}+\tau_{k-1}$ and $X_{-k}=X_{-(k-1)}-\tau_{-k}$, where $\{\tau_k\}_{k\in \Z}$ is an iid sequence of non-negative random variables with bounded second moment.
This point set is interval-stationary, has finite range of dependence but is not stationary up to translation according to Theorem~\ref{th:stat}. We believe there could the following dichotomy in dimension 1: an increment-stationary point set that satisfies Hypothesis~\ref{hypo} for some $\alpha>2$ is either stationary up to translation or interval-stationary up to translation.

\medskip

Before we turn to the proofs of these results, let us focus on a specific class of examples of increment-stationary random point sets given by the image of $\Z^d$ by an ``increment-stationary stochastic diffeomorphism".
This case is of interest for disordered crystals and their thermodynamic limits, and was a motivation for \cite{BLL-07b}.
It will also allow us to emphasize that, although Hypothesis~\ref{hypo} can be interpreted as
a condition on the decay of the correlation between ``stationary increments" (whenever this notion is well-defined), the distance
is the one given by the group action (that is, $k\in \Z^d$ associated with $\theta_k$). In particular Hypothesis~\ref{hypo} 
is \textit{not} a condition on 
the decay of the correlation between the ``stationary increments" wrt to the Euclidean distance.

\subsection{Increment-stationary stochastic diffeomorphisms}

In \cite{BLL-07b}, Blanc, Le Bris and Lions introduced a variant of stochastic homogenization of linear elliptic equations where the diffusion coefficients are random but not necessarily stationary. These diffusion coefficients are obtained using
a stochastic diffeomorphism $\Phi:\R^d\times \Omega \to \R^d$ which is increment-stationary in the following sense: 
For all $k\in \Z^d$, all $x\in \R^d$, and almost every $\omega\in \Omega$,
$$
\nabla \Phi(x+k,\omega)\,=\,\nabla \Phi(x,\theta_k\omega).
$$
Such random fields $\Phi$ allow one to define a specific class of increment-stationary random point sets.
Set $\ell \,:=\,\Phi(\Z^d)$.
Since $\Phi$ is a diffeomorphim, $\ell$ is a simple point process almost surely.
The interest of such a definition is the natural labeling of the points by $\Z^d$.
Define $X_i=\Phi(i)$ for all $i\in \Z^d$ so that $\ell=\cup_{i\in \Z^d}\{X_i\}$.
The increment-stationarity of $\Phi$ then implies that
for all $i,j,k\in \Z^d$ and almost every $\omega \in \Omega$ we have
$$
X_i(\theta_k\omega)-X_j(\theta_k\omega)=X_{i+k}(\omega)-X_{j+k}(\omega),
$$
from which we deduce that $\ell$ is increment-stationary.

\medskip

Such increment-stationary point sets are very specific in the sense that they satisfy the so-called hard-core and non-empty space properties (that is, positive minimal distance between any point $x\in \ell$ and $\ell \setminus \{x\}$, and existence of $R<\infty$ such that any ball of radius $R$
contains at least a point of $\ell$ almost surely) and that they inherit the invariance group of $\Z^d$ (in a statistical way), but not more. In particular such point sets cannot be statistically isotropic (see \cite[Theorem~10]{Alicandro-Cicalese-Gloria-07b}), as opposed e.~g. to the random parking point set (defined in \cite{Penrose-01},
 the properties of which are listed in \cite[Proposition~2.1]{Gloria-Penrose-10}). Anisotropy of the point set is the necessary drawback of the natural labeling of such point sets.
The labeling of $\ell=\Phi(\Z^d)$ has the advantage to allow one to define the notion of increments in the form of the quantities $X_i-X_j$.
As mentioned above, these increments are stationary, and one may consider the associated covariances,
that is, $\cov{(\Phi(k+e_l)-\Phi(k))\cdot e_n}{(\Phi(e_l)-\Phi(0))\cdot e_n}$. 
As we shall show, conditions on the decay of such quantities may ensure the stationarity of the random field $\Phi$, but \textit{not} the stationarity of the point process $\ell=\Phi(\Z^d)$.
More precisely, we have the following counterpart of Proposition~\ref{prop} and Theorem~\ref{th:stat} for increment-stationary random fields.
\begin{prop}\label{prop2}
Let $\Phi:\Z^d\times \Omega\to \R^d$ be an increment-stationary (discrete) random field 
for the group action $\{\theta_k\}_{k\in \Z^d}$, that is, 
such that its discrete gradient $\partial \Phi=(\partial_1 \Phi,\dots,\partial_d \Phi)$, with $\partial_i \Phi:=\Phi(\cdot+e_i)-\Phi(\cdot)$, is stationary: For all $k,z\in \Z^d$ and almost every $\omega\in \Omega$,
$$
\partial \Phi(z+k,\omega)\,=\,\partial \Phi(z,\theta_k\omega).
$$
Assume that $\partial \Phi(0)\in L^2(\Omega, \R^{d\times d})$, and for all $\mu>0$ and $i\in \{1,\dots,d\}$ consider the unique  weak solution $\phi_{\mu,i} \in L^2(\Omega)$ of the equation: For all $\psi\in L^2(\Omega)$,
\begin{equation}\label{eq:corrector}
\expec{\mu \phi_{\mu,i}\psi+D\psi\cdot D\phi_{\mu,i}}\,=\,\expec{D\psi\cdot \partial (\Phi\cdot e_i)(0)}.
\end{equation}
Then, $\Phi$ is stationary up to translation, that is, there exists a unique random vector $\tilde X\in L^2_0(\Omega,\R^d)$ such that $(z,\omega)\mapsto \Phi(z,\omega)+\tilde X(\omega)$ is stationary in the sense that for all $k,z\in \Z^d$ and almost every $\omega\in \Omega$,
$$
\Phi(z+k,\omega)+\tilde X(\omega)\,=\,\Phi(z,\theta_k\omega)+\tilde X(\theta_k\omega),
$$
\emph{if and only if} the family $\{\phi_{\mu,i}\}_{i\in \{1,\dots,d\}}$ is bounded in $L^2(\Omega)$ uniformly wrt $\mu>0$.
\end{prop}
The conditions corresponding to Hypothesis~\ref{hypo} are now
\begin{hypo}[Decay of order $\alpha>0$]\label{hypo2}
The random field $\Phi:\Z^d\times \Omega\to \R^d$ is increment-stationary for the ergodic group action $\{\theta_k\}_{k\in \Z^d}$,
and there exists $\alpha>0$ such that for all $k\in \Z^d$ and $l,l',n\in \{1,\dots,d\}$
$$
\cov{\partial_l \Phi(k)\cdot e_n}{\partial_{l'}\Phi(0)\cdot e_n} \,\lesssim \,\frac{1}{1+|k|^\alpha},
$$
where $\lesssim$ means $\leq$ up to a universal multiplicative constant.
\end{hypo}
\begin{theo}\label{th:stat2} We have:
\begin{itemize}
\item For $d=1$ and $d=2$: There exists an increment-stationary random field $\Phi:\Z^d\times \Omega\to \R^d$ satisfying Hypothesis~\ref{hypo2} with finite range of dependence and which is not stationary up to translation.
\item For all $d>2$: If $\Phi:\Z^d\times \Omega\to \R^d$ is increment-stationary and satisfies Hypothesis~\ref{hypo2} for some $\alpha>2$, then it is stationary up to translation.
\end{itemize}
\end{theo}
As a direct corollary of Theorem~\ref{th:stat2}, we have the following result on random Lipschitz fields with stationary gradients:
\begin{coro}\label{th:diffeo}
Let $d>2$ and $\Phi:\R^d\times \Omega \to \R^d$ be a random Lipschitz field such that its (continuum) gradient $\nabla \Phi$ is stationary and uniformly bounded.
If there exists $\alpha>2$ such that for all $k\in \Z^d,x\in [0,1)^d$ and $l,l',n\in \{1,\dots,d\}$
$$
\cov{\partial_l\Phi(x+k)\cdot e_n}{\partial_{l'}\Phi(x)\cdot e_n}\,\lesssim \, \frac{1}{1+|k|^\alpha},
$$
then there exists a $[0,1)^d$-periodic random field $\tilde X\in W^{1,\infty}_\per([0,1)^d,L^2_0(\Omega,\R^d))$ such that $(x,\omega)\mapsto \Phi(x,\omega)+\tilde X(x,\omega)$ is a stationary field: For all $x\in \R^d$, $k\in \Z^d$, and almost every $\omega\in \Omega$,
$$
\Phi(x+k,\omega)+\tilde X(x+k,\omega)\,=\,\Phi(x,\theta_k\omega)+\tilde X(x,\theta_k \omega).
$$
\end{coro}
To conclude, let us compare the two notions of stationarity for a random field $\Phi$
 and for the associated random point set $\ell=\Phi(\Z^d)$.
The following result shows that these notions are essentially incompatible.
\begin{prop}\label{prop:incomp}
Let $\Phi:\Z^d\times \Omega \to \R^d$ be an increment-stationary (discrete) injective random field 
for the ergodic group action $\{\theta_k\}_{k\in \Z^d}$, that is, 
such that its discrete gradient $\partial \Phi=(\partial_1 \Phi,\dots,\partial_d \Phi)$, with $\partial_i \Phi:=\Phi(\cdot+e_i)-\Phi(\cdot)$, is stationary.
If $\Phi$ and $\ell=\Phi(\Z^d)$ are stationary up to translation, then $\Phi$ is linear and $\ell$ is periodic up to a random translation. Indeed, we then have 
$\Phi(\cdot,\omega):x\mapsto \expec{\partial \Phi}x+\Phi(0,\omega)$, and
$\ell(\omega)=\expec{\partial \Phi}\Z^d+\Phi(0,\omega)$ for almost every $\omega\in \Omega$.
\end{prop}
This proposition illustrates that even in the case when one may properly define the notion of ``stationary increments",
 Hypothesis~\ref{hypo} cannot be turned into a condition on the decay of the correlation between the ``stationary increments" wrt to the Euclidean distance (which corresponds to Hypotheses~\ref{hypo2}).

\section{Proofs}

The proofs of Proposition~\ref{prop2} and Theorem~\ref{th:stat2} are straightforward adaptations of the proofs of  Proposition~\ref{prop} and Theorem~\ref{th:stat}, and we only prove the latter.

\subsection{Proof of Proposition~\ref{prop}}

We split the proof into two steps.

\medskip

\step{1} Proof that boundedness of $\{\phi_{\mu,i}\}_{i\in \{1,\dots,d\}}$ implies stationarity up to translation.

\noindent Let $\ell$ be an increment-stationary random point set.
Since the action group is ergodic, it generates a Weyl decomposition of $L^2(\Omega,\R^d)$ into potential fields (that is, the closure in $L^2(\Omega,\R^d)$ of gradient fields) and solenoidal fields (that is, the vector fields that are orthogonal to potential fields for the $L^2(\Omega,\R^d)$-scalar product), see for instance \cite[Lemma~7.3]{JKO-94}, the adaptation of which is straightforward in the case of a discrete group action.
Taking $\psi=\phi_{\mu,i}$ in \eqref{eq:corrector} yields 
$$
\mu\expec{\phi_{\mu,i}^2}+\expec{|D\phi_{\mu,i}|^2} \,=\,\expec{D\phi_{\mu,i}\cdot \zeta_i},
$$
which turns, by Cauchy-Schwarz' inequality and the assumption $\zeta_i\in L^2(\Omega,\R^d)$, into the energy estimate 
$$
\mu\expec{\phi_{\mu,i}^2}+\expec{|D\phi_{\mu,i}|^2}\lesssim 1.
$$
By the Banach-Alaoglu theorem the sequence $D\phi_{\mu,i}$ is weakly compact, and converges weakly in $L^2(\Omega,\R^d)$ to some potential field $\chi_i \in L^2(\Omega,\R^d)$ up to extraction.
Passing to the limit along the subsequence in the defining equation \eqref{eq:corrector} for $\phi_{\mu,i}$ and using the a priori estimate $\mu\expec{\phi_{\mu,i}^2}\lesssim 1$ yield
for all $\psi\in L^2(\Omega)$,
$$
\expec{D\psi\cdot \chi_i}\,=\, \expec{D\psi \cdot \zeta_i}.
$$
The above identity for arbitrary $\psi\in L^2(\Omega)$ shows that the $L^2(\Omega)$-projections of $\chi_i$ and $\zeta_i$ onto potential fields coincide.
Since $\zeta_i = ((Z_{e_1}-Z_0)\cdot e_i,\dots,(Z_{e_d}-Z_0)\cdot e_i)=D(Z\cdot e_i)$ and $\chi_i$ (as limit of potential fields $D\phi_{\mu,i}$) are potential fields themselves, this implies  $\chi_i=\zeta_i$
and yields the uniqueness of the limit.
Assume in addition that the family $\phi_{\mu,i}$ is bounded
in $L^2(\Omega)$ uniformly wrt $\mu$. 
Then, by weak compactness, there exists some $\phi_i\in L^2(\Omega)$ such that, up to extraction, $\phi_{\mu,i}$ converges to $\phi_i$ weakly in $L^2(\Omega)$. Note that $\expec{\phi_i}=0$ since
$\expec{\phi_{\mu,i}}=0$ for all $\mu>0$.
Combined with the argument above, this shows that $D\phi_i=\zeta_i$, and this implies in turn the uniqueness of $\phi_i$ and the convergence of the entire sequence.
Indeed, let $\varphi\in L^2(\Omega)$ be such that $\expec{\varphi}=\expec{\phi_i}=0$ and $D\varphi =\zeta_i$. Then $\varphi-\phi_i\in L^2(\Omega)$ is such that
$\expec{|D(\varphi-\phi_i)|^2}=0$. This implies by ergodicity that $\varphi-\phi_i$ is constant, and therefore $\varphi=\phi_i$ by the mean-free condition.

We then define $\tilde Y:\Omega \to \R^d, \omega \mapsto \sum_{i=1}^d \phi_i(\omega)e_i$,
and for all $i\in \{1,\dots,d\}$, 
we set $T_i:=\expec{Y_{e_i}}$.
It remains to check that $\ell+\tilde Y$ is stationary.

Using that $Y_0\equiv 0$ one can decompose $Y_k$ as a sum of differences along the $d$ canonical directions, i.~e.
\begin{multline*}
Y_k\,=\, \sum_{i_1=1}^{k_1} \big( Y_{k_1+1-i_1,k_2,\dots,k_d}-Y_{k_1-i_1,k_2,\dots,k_d} \big)+\sum_{i_2=1}^{k_2} \big(Y_{0,k_2+1-i_2,k_3,\dots,k_d}-Y_{0,k_2-i_2,k_3,\dots,k_d}\big)
\\
+\cdots +\sum_{i_d=1}^{k_d} \big(Y_{0,\dots,0,k_d+1-i_d}-Y_{0,\dots,0,k_d-i_d}\big).
\end{multline*}
By stationarity of the increments and definition of $\{Z_{e_i}\}_{i\in \{1,\dots,d\}}$, this yields 
\begin{eqnarray*}
\expec{Y_k}&=& \sum_{i=1}^d k_i\expec{Y_{e_i}},
\\
Y_k(\omega)-\expec{Y_k}&=& \sum_{i_1=1}^{k_1} Z_{e_1}(\theta_{k_1-i_1,k_2,\dots,k_d}\omega) +\sum_{i_2=1}^{k_2} Z_{e_2}(\theta_{0,k_2-i_2,k_3,\dots,k_d}\omega)\\
&&+\cdots+\sum_{i_d=1}^{k_d} Z_{e_d}(\theta_{0,\dots,0,k_d-i_d}\omega).
\end{eqnarray*}
Using $Z_0\equiv 0$ and $D\phi_i=\zeta_i$, this turns into
\begin{multline*}
Y_k(\omega)-\expec{Y_k}\,=\, \sum_{i_1=1}^{k_1} \sum_{l=1}^d D_1 \phi_l (\theta_{k_1-i_1,k_2,\dots,k_d}\omega)e_l +\sum_{i_2=1}^{k_2} \sum_{l=1}^dD_2 \phi_l(\theta_{0,k_2-i_2,k_3,\dots,k_d}\omega)e_l
\\
+\cdots+\sum_{i_d=1}^{k_d} \sum_{l=1}^d D_d \phi_l(\theta_{0,\dots,0,k_d-i_d}\omega)e_d.
\end{multline*}
By definition of the difference operators $D_i$, terms cancel two by two, and the sum simplifies to 
\begin{equation*}
Y_k(\omega)-\expec{Y_k}\,=\, \sum_{l=1}^d \phi_l(\theta_{k}\omega)e_l - \sum_{l=1}^d \phi_l(\omega)e_l.
\end{equation*}
This implies the desired property by the choice of $\tilde Y$ and $T_i$: For all $k\in \Z^d$,
$$
\ell(\theta_k\omega)+\tilde Y(\theta_k\omega)\,=\,\ell(\omega)+\tilde Y(\omega)+\sum_{i=1}^d k_i T_i.
$$

\medskip

\step{2} Proof that stationarity up to translation implies boundedness of $\{\phi_{\mu,i}\}_{i\in \{1,\dots,d\}}$. 

\noindent On the one hand, since $\ell$ is stationary up to translation, there exist $\tilde Y\in L^2(\Omega,\R^d)$ and $\{T_i\}_{1\leq i\leq d}\in \R^d$
such that for all $k\in \Z^d$ and almost all $\omega$, 
$$
\ell(\theta_k\omega)+\tilde Y(\theta_k\omega)=\ell(\omega)+\tilde Y(\omega)+\sum_{i=1}^d k_i T_i.
$$
On the other hand, increment-stationarity implies there exists $Y_k\in L^2(\Omega,\R^d)$ such that $\ell(\theta_k \omega)=\ell(\omega)+Y_k(\omega)$. Since $\ell$ is assumed to be non-degenerate, this implies 
that $Y_k(\omega)=\tilde Y(\omega)-\tilde Y(\theta_k\omega)+\sum_{i=1}^d k_i T_i$.
Hence, $\expec{Y_k}=\sum_{i=1}^d k_i T_i$ and $Z_{e_l}$ takes the form
$$
Z_{e_l}(\omega)\,=\,\tilde Y(\theta_{e_l}\omega)-\tilde Y(\omega).
$$
Recall that for all $\mu>0$ and $i\in \{1,\dots, d\}$, $\phi_{\mu,i}\in L^2(\Omega)$ is solution of: For all $\psi\in L^2(\Omega)$,
$$
\mu \expec{\phi_{\mu,i}\psi}+\expec{D\psi\cdot D\phi_{\mu,i}}\,=\,\expec{D\psi\cdot \zeta_i}
$$
with $\zeta_i=(Z_{e_1}\cdot e_i, \dots, Z_{e_d}\cdot e_i)$.
For all $i\in \{1,\dots, d\}$,  set $\tilde \phi_i:=\tilde Y \cdot e_i \in L^2(\Omega)$, so that $D\tilde \phi_i=\zeta_i$.
We then have for all $\psi\in L^2(\Omega)$,
$$
\mu\expec{(\phi_{\mu,i}-\tilde \phi_i)\psi}+\expec{D\psi\cdot D(\phi_{\mu,i}-\tilde \phi_i)}\,=\,-\mu\expec{\tilde \phi_{i}\psi},
$$
whence the a priori estimate 
$$
\expec{(\phi_{\mu,i}-\tilde \phi_i)^2}\,\leq \, \expec{\tilde \phi_i^2}=\expec{(\tilde Y\cdot e_i)^2},
$$
which yields the desired uniform-in-$\mu$ boundedness estimate by the triangle inequality:
$$
\expec{\phi_{\mu,i}^2}\,\lesssim \, \expec{(\tilde Y\cdot e_i)^2}.
$$

\subsection{Proof of Theorem~\ref{th:stat}}

Let $\ell$ be an increment-stationary random point set, and for all $i\in \{1,\dots,d\}$ and $\mu>0$, let $\zeta_i$ and
$\phi_{\mu,i}$ be as in Proposition~\ref{prop}.
We first derive an integral representation for $\phi_{\mu,i}$ in physical space, then treat the case 
$d\leq 2$ in Step~2, and the case $d>2$ in Step~3.

\medskip

\step{1} Green representation formula for $\phi_{\mu,i}$.

\noindent 
In this step we derive a Green representation formula for $\phi_{\mu,i}$, see \cite[Lemma~2.6]{Gloria-Otto-09}.
Equation~\eqref{eq:corrector} indeed admits an equivalent form in the physical space.
Let $\bar \phi_{\mu,i},\bar \zeta_i:\Z^d\times \Omega\to \R$ be the stationary extensions of $\phi_{\mu,i}$ and $\zeta_i$, that is, 
the random fields defined by $\bar \phi_{\mu,i}(k,\omega):=\phi_{\mu,i}(\theta_k \omega)$ and $\bar \zeta_{i}(k,\omega):=\zeta_{i}(\theta_k \omega)$, respectively. 
Then, $\bar \phi_{\mu,i}$ solves almost surely the elliptic PDE
\begin{equation}\label{eq:corr-space}
\mu \bar \phi_{\mu,i} -\triangle \bar \phi_{\mu,i} \,=\, \partial^* \cdot \bar \zeta_i \qquad \text{ in }\Z^d,
\end{equation}
where $\partial$ is the forward discrete gradient, $\partial^*\cdot$ the backward discrete divergence, and $-\triangle\,:=\,-\partial^*\cdot \partial$ the discrete
Laplace operator on $\Z^d$.

Let $G_\mu:\Z^d\to \R$ denote the Green's function associated with the elliptic operator $\mu-\triangle$ on $\Z^d$, that is, the only solution in $L^2(\Z^d)$ of
$$
\mu G_\mu(x)-\triangle G_\mu(x)\,=\,\delta(x),
$$
where $\delta$ is such that $\delta(0)=1$ and $\delta(x)=0$ for all $x\neq 0$ (the existence and uniqueness of $G_\mu$ follows from the Riesz representation theorem).

Testing equation \eqref{eq:corr-space} with $y\mapsto G_\mu(y-x)$ yields the desired Green representation formula
$$
\bar \phi_{\mu,i}(x)\,=\,\int_{\Z^d} \partial G_\mu(y-x) \cdot \bar \zeta_i(y)dy,
$$
where $\int_{\Z^d}dy$ stands for the sum over $y\in \Z^d$.

\medskip

\step{2} Case $d\leq 2$.

\noindent In this step we prove that even in the case when $\{\bar \zeta_i(z)\}_{z\in \Z^d}$ is a field of independent and identically distributed (iid) variables, the family $\expec{\phi_{\mu,i}^2}$ may be unbounded in $\mu$.
Consider in particular the field $Y_k$ characterized by: $Y_0\equiv 0$ and $Y_{e_l}\circ \theta_k=Y_{k+e_l}-Y_{k}=a_l(k)e_l$, where $\{a_l(k)\}_{l\in \{1,\dots,d\},k\in \Z^d}$ are iid variables following the law of some $a\in L^2(\Omega)$. 
Then, the random point set $\ell$ satisfies Hypothesis~\ref{hypo} for any $\alpha>0$ (it has finite correlation-length).
In view of Step~1, we have
\begin{equation}\label{eq:2nd-moment}
\expec{\phi_{\mu,i}^2}\,=\,\expec{(\bar \phi_{\mu,i}(0))^2}\,=\,\int_{\Z^d}\int_{\Z^d} \partial G_\mu(y) \otimes \partial G_\mu(y') : \expec{\bar \zeta_i(y)\otimes  \bar \zeta_i(y')}dydy'.
\end{equation}
Since $\bar \zeta_i(y)=(a_i(y)-\expec{a})e_i$, the sum reduces by independence to 
\begin{equation}\label{eq:phimu-d=2}
\expec{\phi_{\mu,i}^2}\,=\,\var{a} \int_{\Z^d} (\partial_i G_\mu(y))^2 dy.
\end{equation}
Since $\partial G_\mu$ converges locally to the gradient $\partial G$ of the whole space Green's function of the discrete Laplace operator $-\triangle$ on $\Z^d$ as $\mu\downarrow 0$,
the RHS of \eqref{eq:phimu-d=2} cannot be bounded in dimensions $d\leq 2$. If it were, this would imply that $\partial G\in L^2(\Z^d)$, which is not true for $d\leq 2$ (as can be seen in Fourier space \cite{Martinsson-02},
or by comparison to the large scale behavior of the continuum Green's function).
This qualitative behavior is enough for the proof of Theorem~\ref{th:stat}. To be more quantitative, one indeed expects $\expec{\phi_{\mu,i}^2}\sim \mu^{-\frac{1}{2}}$ for $d=1$, and $\expec{\phi_{\mu,i}^2}\sim |\ln \mu|$ for $d=2$. The proof of these estimates would require
a more careful analysis of the Green's functions.

\medskip

\step{3} Case $d>2$.

\noindent The starting point in dimensions $d>2$ is again \eqref{eq:2nd-moment} which in view of Hypothesis~\ref{hypo} implies
$$
\expec{\phi_{\mu,i}^2}\,\lesssim \, \int_{\Z^d}\int_{\Z^d} | \partial G_\mu(y)|| \partial G_\mu(y')| \frac{1}{1+|y-y'|^\alpha}dydy'.
$$
Without loss of generality, we assume in addition that $\alpha<d$.
We shall use the following uniform-in-$\mu$ bounds on $\partial G_\mu$: $\|\partial G_\mu\|_{L^\infty(\Z^d)} \lesssim 1$ (cf. \cite[Corollary~2.3]{Gloria-Otto-09}), and for all exponents $1\leq  p<\infty$, all $i\in \N$ and $d>2$,
\begin{equation}\label{eq:dyad}
\int_{2^i <|y|\leq 2^{i+1}} |\partial G_\mu(y)|^p dy \,\lesssim\, (2^i)^d (2^i)^{p(1-d)},
\end{equation}
which is optimal in terms of scaling.
This estimate is standard and relies on the $L^p$-regularity theory for the operator $\mu-\triangle$.
For a proof, we refer to \cite[Lemma~2.9]{Gloria-Otto-09}, which treats in addition the variable-coefficients case using the perturbation approach by Meyers.
Indeed, Steps~3--6 of that proof show that if for some $p>2$ the operator has an $L^p$-regularity theory, then \eqref{eq:dyad} holds, whereas Step~1 shows that
$\mu-\triangle$ has an $L^p$-regularity theory for all $1<p<\infty$. The case $1\leq p< 2$ in \eqref{eq:dyad} follows from the estimate for $p=2$ by H\"older's inequality.

We now prove the boundedness of $\expec{\phi_{\mu,i}^2}$ if $\alpha>2$ by estimating the integrals using a doubly dyadic decomposition of $\Z^d\times \Z^d$.
Note that the exponent $\alpha=2$ is borderline in terms of integrability.

We write the integral as:
\begin{subequations}
\begin{eqnarray}
\lefteqn{\int_{\Z^d}\int_{\Z^d} | \partial G_\mu(y)|| \partial G_\mu(y')| \frac{1}{1+|y-y'|^\alpha}dydy' } \nonumber \\
&=& \int_{|y|\leq 2}| \partial G_\mu(y)|  \int_{|y-y'|\leq 2}| \partial G_\mu(y')| \frac{1}{1+|y-y'|^\alpha}dydy' \label{eq:rhs1}\\
&&+\int_{|y|\leq 2}| \partial G_\mu(y)|\sum_{j=1}^{\infty} \int_{2^j <|y-y'|\leq 2^{j+1}}| \partial G_\mu(y')| \frac{1}{1+|y-y'|^\alpha}dydy' \label{eq:rhs2}\\
&&+\sum_{i=1}^\infty \int_{2^i <|y|\leq 2^{i+1}} | \partial G_\mu(y)| \int_{|y-y'|\leq 2}| \partial G_\mu(y')| \frac{1}{1+|y-y'|^\alpha}dydy'  \label{eq:rhs3}\\
&&+\sum_{i=1}^\infty \int_{2^i <|y|\leq 2^{i+1}} | \partial G_\mu(y)| \sum_{j=1}^{\infty} \int_{2^j <|y-y'|\leq 2^{j+1}}| \partial G_\mu(y')| \frac{1}{1+|y-y'|^\alpha}dydy' \label{eq:rhs4}.
\end{eqnarray}
\end{subequations}
By the uniform bound on $\|\partial G_\mu\|_{L^\infty(\Z^d)}$, the RHS term \eqref{eq:rhs1} is of order $1$.
For the RHS term \eqref{eq:rhs2}, this uniform bound and the triangle inequality yield
$$
\eqref{eq:rhs2}\,\lesssim\,1+\sum_{j=1}^{\infty} \int_{2^j-2 <|y'|\leq 2^{j+1}+2}| \partial G_\mu(y')| \frac{1}{1+|y'|^\alpha}dy'.
$$
Estimate~\eqref{eq:dyad} for $p=1$ then yields that $\eqref{eq:rhs2}\,\lesssim\, 1$ since $\alpha>1$.
The proof of the boundedness of the RHS term \eqref{eq:rhs3} is similar.
The most subtle part is the RHS term \eqref{eq:rhs4}. We split the double sum into two parts: $\sum_{i=1}^\infty \sum_{j\leq i}$ and $\sum_{i=1}^\infty\sum_{j>i}$, and start with the latter.
If $j> i$, and $y,y'$ are such that $2^i <|y|\leq 2^{i+1}$ and  $2^j <|y-y'|\leq 2^{j+1}$, then $2^{j-1}<|y'|\leq 2^{j+2}$.
Hence, by \eqref{eq:dyad} for $q=1$,
\begin{eqnarray*}
\lefteqn{\int_{2^i <|y|\leq 2^{i+1}} | \partial G_\mu(y)| \int_{2^j <|y-y'|\leq 2^{j+1}}| \partial G_\mu(y')| \frac{1}{1+|y-y'|^\alpha}dydy' }
\\
&\leq& (2^{j})^{-\alpha}\int_{2^i <|y|\leq 2^{i+1}} | \partial G_\mu(y)|dy \int_{2^{j-1} <|y'|\leq 2^{j+2}}| \partial G_\mu(y')|dy' \\
\\
&\lesssim&  (2^i)^{d-(d-1)} (2^j)^{d-(d-1)-\alpha} \,=\,  (2^i)  (2^j)^{1-\alpha}.
\end{eqnarray*}
Since $\alpha>1$, summing over $j>i$ yields
\begin{equation*}
\int_{2^i <|y|\leq 2^{i+1}} | \partial G_\mu(y)| \sum_{j> i} \int_{2^j <|y-y'|\leq 2^{j+1}}| \partial G_\mu(y')| \frac{1}{1+|y-y'|^\alpha}dydy' \,\lesssim \, (2^i)^{2-\alpha},
\end{equation*}
and therefore, using that $\alpha>2$,
\begin{equation}\label{eq:part2}
\sum_{i=1}^\infty \int_{2^i <|y|\leq 2^{i+1}} | \partial G_\mu(y)| \sum_{j> i} \int_{2^j <|y-y'|\leq 2^{j+1}}| \partial G_\mu(y')| \frac{1}{1+|y-y'|^\alpha}dydy'  \,\lesssim\, 1.
\end{equation}
We now treat the sum $\sum_{i=1}^\infty \sum_{j\leq i}$.
If $j\leq i$, and $y,y'$ are such that $2^i <|y|\leq 2^{i+1}$ and  $2^j <|y-y'|\leq 2^{j+1}$, then $2^{i-1}<|y'|\leq 2^{i+2}$.
Let $q>1$ be such that that $\frac{d}{q}>\alpha$ (which is possible since $d>\alpha$), and $p>1$ be the associated dual exponent, i.~e. $\frac{1}{p}+\frac{1}{q}=1$. 
Then, by H\"older's inequality with exponents $(p,q)$, and \eqref{eq:dyad} with exponents 1 and $p$, we have
\begin{eqnarray*}
\lefteqn{\int_{2^i <|y|\leq 2^{i+1}} | \partial G_\mu(y)| \int_{2^j <|y-y'|\leq 2^{j+1}}| \partial G_\mu(y')| \frac{1}{1+|y-y'|^\alpha}dydy' }
\\
&\leq& \int_{2^i <|y|\leq 2^{i+1}} | \partial G_\mu(y)|dy \left(\int_{2^{i-1} <|y'|\leq 2^{i+2}}| \partial G_\mu(y')|^pdy'\right)^{\frac{1}{p}} \\
&& \qquad \times
\left(\int_{2^j <|y-y'|\leq 2^{j+1}} \frac{1}{1+|y-y'|^{q\alpha}}dy'\right)^{\frac{1}{q}} 
\\
&\lesssim&  (2^i)^{d-(d-1)} (2^i)^{\frac{1}{p}(d-p(d-1))} (2^j)^{\frac{1}{q}(d-q\alpha)} \,=\, (2^i)^{2-d(1-\frac{1}{p})}  (2^j)^{\frac{d}{q}-\alpha}.
\end{eqnarray*}
Summing over $j\leq i$ and using that $\frac{d}{q}-\alpha>0$ then yields
\begin{multline*}
{\int_{2^i <|y|\leq 2^{i+1}} | \partial G_\mu(y)| \sum_{j\leq i} \int_{2^j <|y-y'|\leq 2^{j+1}}| \partial G_\mu(y')| \frac{1}{1+|y-y'|^\alpha}dydy' }
\\
\lesssim\,   (2^i)^{2-d(1-\frac{1}{p})} \sum_{j\leq i}  (2^j)^{\frac{d}{q}-\alpha} \,
\lesssim \, (2^i)^{2-d(1-\frac{1}{p})+\frac{d}{q}-\alpha}\,=\,(2^i)^{2-\alpha}. 
\end{multline*}
Since $\alpha>2$, this yields a bound for the first sum $\sum_{i=1}^\infty \sum_{j\leq i}$:
\begin{equation}\label{eq:part1}
\sum_{i=1}^\infty \int_{2^i <|y|\leq 2^{i+1}} | \partial G_\mu(y)| \sum_{j\leq i} \int_{2^j <|y-y'|\leq 2^{j+1}}| \partial G_\mu(y')| \frac{1}{1+|y-y'|^\alpha}dydy'  \,\lesssim\, 1.
\end{equation}
The combination of \eqref{eq:part2} and \eqref{eq:part1} shows that $\eqref{eq:rhs4}\lesssim 1$, which, combined
with the estimates of \eqref{eq:rhs1}, \eqref{eq:rhs2}, and \eqref{eq:rhs3}, implies  the uniform-in-$\mu$ bound
$$
\int_{\Z^d}\int_{\Z^d} | \partial G_\mu(y)|| \partial G_\mu(y')| \frac{1}{1+|y-y'|^\alpha}dydy' \,\lesssim\, 1,
$$
as desired.

\subsection{Proof of Corollary~\ref{th:diffeo}}

For all $x\in \R^d$, consider the random field $\Phi_x:\Z^d\times \Omega : (k,\omega)\mapsto \Phi(x+k,\omega)$.
By Theorem~\ref{th:stat2}, there exists some random vector $\tilde{X}(x,\cdot)\in L^2_0(\Omega)$ such that $(z,\omega)\mapsto \Phi_x(z,\omega)+\tilde X(x,\omega)$ is stationary: For all $z,k\in \Z^d$ and almost every $\omega\in \Omega$,
$$
\Phi_x(z+k,\omega)+\tilde X(x,\omega)\,=\,\Phi_x(z,\theta_k\omega)+\tilde X(x,\theta_k\omega),
$$
which we may rewrite by definition of $\Phi_x$ as 
\begin{equation}\label{eq:sto-diff-stat}
\Phi(x+z+k,\omega)+\tilde X(x,\omega)\,=\,\Phi(x+z,\theta_k\omega)+\tilde X(x,\theta_k\omega).
\end{equation}
Since $\Phi_{x+k}(k')=\Phi_{x}(k+k')$ for all $k,k'\in \Z^d$, the uniqueness of $\tilde X$ (cf. uniqueness of $\phi_i$ in Step~1 of the proof of Proposition~\ref{prop}) shows that 
$\tilde X(x+k,\cdot)\,=\,\tilde X(x,\cdot)$ for all $k\in \Z^d$ and $x\in \R^d$.
Hence, \eqref{eq:sto-diff-stat} turns into: For all $x\in \R^d$, there exists a set of full measure $\Omega_x$ such that
for all $z,k\in \Z^d$ and $\omega\in \Omega_x$, we have
\begin{equation*}
\Phi(x+z+k,\omega)+\tilde X(x+z+k,\omega)\,=\,\Phi(x+z,\theta_k\omega)+\tilde X(x+z,\theta_k\omega).
\end{equation*}
To conclude, it remains to prove the measurability of $\tilde X:\R^d\times \Omega\to \R^d$ (where $\R^d$ is endowed
with the Borel sets).
It is enough to show that $\tilde X$ is a Carath\'eodory function, since Carath\'eodory functions are equivalent to Borel functions. Recall that $\tilde X$ is a Carath\'eodory function if for almost all $\omega\in \Omega$, $x\mapsto \tilde X(x,\omega)$ is continuous, and if for all $x\in \R^d$, $\omega\mapsto \tilde X(x,\omega)$ is measurable.
The measurability of $\tilde X(x,\cdot)$ follows form the definition of $\tilde X$. It remains to prove the continuity, which we do in the form of a Lipschitz estimate.
There exists a set $\Omega'\in \calF$ of full measure such that for all $x\in \Q^d$, $z,k\in \Z^d$ and $\omega\in \Omega'$,
\begin{equation}\label{eq:sto-diff-stat2}
\Phi(x+z,\omega)+\tilde{X}(x+z,\omega)\,=\,\Phi(x,\theta_k\omega)+\tilde{X}(x,\theta_k\omega).
\end{equation}
The uniform Lipschitz assumption on $\Phi$ yields: There exists $C<\infty$ such that for all $k\in \Z^d$ and $x,h\in \R^d$, $|\Phi(x+h+k,\cdot)-\Phi(x+k,\cdot)|\leq C |h|$. 
Hence, substracting \eqref{eq:sto-diff-stat2} once with $x\leadsto x+h$ and once with $x$ implies that for all $x,h \in \Q^d$,  $k\in \Z^d$, and $\omega\in \Omega'$, 
$$
\tilde{X}(x+h,\omega)-\tilde{X}(x,\omega) - \Big(  \tilde{X}(x+h,\theta_k \omega)-\tilde{X}(x,\theta_k\omega)\Big) \,\leq \, C|h|.
$$
By summation over $k$, this yields for all $N\in \N$
$$
\tilde X(x+h,\omega)-\tilde X(x,\omega)-\frac{1}{\#([-N,N)\cap \Z)^d} \sum_{k\in ([-N,N)\cap \Z)^d}\tilde X(x+h,\theta_k \omega)-\tilde X(x,\theta_k\omega) \,\leq\, C|h|.
$$
By the ergodic theorem, and since $\expec{\tilde X(x,\cdot)}=\expec{\tilde X(x+h,\cdot)}=0$, there exists some $\Omega''\in \calF$ with full measure such that for all $x \in \Q^d$, $h\in \Q^d$, $k\in \Z^d$ and $\omega\in \Omega''$ the limit $N\uparrow \infty$ yields
$$
\tilde X(x+h,\omega)-\tilde X(x,\omega) \,\leq\, C|h|.
$$
By symmetry, this implies
$$
|\tilde X(x+h,\omega)-\tilde X(x,\omega)|\,\leq\, C|h|,
$$
so that $\tilde X|_{\Q^d\times \Omega''}$ can be extended to a Lipschitz function on $\R^d$ for all $\omega\in \Omega''$. 
Hence, $\tilde X$ is a Carath\'eodory function, as desired.

\subsection{Proof of Proposition~\ref{prop:incomp}}

We split the proof into two steps, using the stationarity of $\ell+\tilde Y=\Phi(\Z^d)+\tilde Y$, and of $\Phi+\tilde X$, respectively.
Recall that $T=\expec{\partial \Phi}$.

\medskip

\step{1} Stationarity of $\ell$.

\noindent Since $\ell$ is stationary up to translation, there exists some $\tilde Y \in L^2_0(\Omega)$ such that for all $k\in \Z^d$ and almost every $\omega \in \Omega$,
$$
\ell(\theta_k\omega)+\tilde Y(\theta_k\omega)\,=\,\ell(\omega)+\tilde Y(\omega)+Tk.
$$
Every point of $\ell(\theta_k\omega)+\tilde Y(\theta_k\omega)$ is mapped to a point of $\ell(\omega)+\tilde Y(\omega)+Tk$. Hence for almost all $\omega \in \Omega$, there exists a function $\gamma_\omega:\Z^d\times \Z^d\to \Z^d$
such that for all $k\in \Z^d$, $\gamma_\omega(k,\cdot)$ is bijective on $\Z^d$ and satisfies
$$
X_y(\theta_k \omega)+\tilde Y(\theta_k\omega)\,=\,X_{\gamma(k,\omega,y)}(\omega)+\tilde Y(\omega)+Tk,
$$
where $X_y(\omega)=\Phi(y,\omega)$ as before.
Combined with the stationarity of the increment $X_y-X_0$, this turns into
\begin{eqnarray}
X_y(\theta_k\omega)+\tilde Y(\theta_k\omega)&=&X_0(\theta_k\omega)+\tilde Y(\theta_k\omega)+X_y(\theta_k\omega)-X_0(\theta_k\omega)\nonumber \\
&=&X_y(\omega)+\tilde Y(\omega)+Tk+X_{\gamma(k,\omega,0)}(\omega)-X_0(\omega)\label{eq:station-ell-1}.
\end{eqnarray}

\medskip

\step{2} Stationarity of $\Phi$ and conclusion.

\noindent Since $\Phi$ is stationary up to translation, there exists $\tilde X\in L^2_0(\Omega)$ such that for all $y,k\in \Z^d$ and almost every $\omega\in \Omega$,
\begin{eqnarray*}
X_{y+k}(\omega)+\tilde X(\omega)&=&X_y(\theta_k\omega)+\tilde X(\theta_k\omega).
\end{eqnarray*}
Combined with \eqref{eq:station-ell-1} this yields
\begin{equation*}
X_{y+k}(\omega)-X_y(\omega)\,=\,Tk+\tilde X(\theta_k\omega)-\tilde X(\omega)+\tilde Y(\omega)-\tilde Y(\theta_k \omega)+X_{\gamma(k,\omega,0)}(\omega)-X_0(\omega).
\end{equation*}
Since the RHS does not depend on $y$, the increment $X_{y+k}-X_y$ does not depend on $y$ either.
By the ergodic theorem, and since the increment is stationary, this yields
for almost every $\omega \in \Omega$ and for all $y,k\in \Z^d$,
\begin{multline*}
X_{y+k}-X_y \,=\, \lim_{R\to \infty} \frac{1}{\# ([-R,R)\cap \Z)^d} \sum_{z\in ([-R,R)\cap \Z)^d} X_{z+k}-X_z \\
\,=\,\expec{X_{k}-X_0}
\,=\,\expec{\partial \Phi}k\,=\,Tk.
\end{multline*}
Hence, for almost every $\omega \in \Omega$ and all $y \in \Z^d$,
$$
X_y(\omega)\,=\,X_0(\omega)+Ty,
$$
so that we obtain
\begin{eqnarray*}
\Phi(y,\omega)&=&\expec{\partial \Phi}y+\Phi(0,\omega),\\
\ell(\omega)&=&\expec{\partial \Phi}\Z^d+\Phi(0,\omega),
\end{eqnarray*}
as desired.

\section*{Acknowledgements}
The author is grateful to Xavier Blanc for his constructive comments
on a early version of this paper and to the anonymous referee for insightful suggestions to improve readability.
The author acknowledges financial support from ANR-10-JCJC-0106 AMAM and from the European Research Council under
the European Community's Seventh Framework Programme (FP7/2014-2019 Grant Agreement
QUANTHOM 335410).

\bibliographystyle{plain}

\end{document}